\documentclass{amsart}

\usepackage{amssymb,latexsym}

\newtheorem{theorem}{Theorem}
\newcommand{\bt}{\begin{theorem}}
\newcommand{\et}{\end{theorem}}
\newtheorem{lemma}{Lemma}
\newcommand{\bl}{\begin{lemma}}
\newcommand{\el}{\end{lemma}}
\newtheorem{corollary}{Corollary}
\newcommand{\bc}{\begin{corollary}}
\newcommand{\ec}{\end{corollary}}
\newtheorem{problem}{Problem}
\newcommand{\bprob}{\begin{problem}}
\newcommand{\eprob}{\end{problem}}
\newcommand{\beq}{\begin{equation}}
\newcommand{\eeq}{\end{equation}}
\newcommand{\benum}{\begin{enumerate}}
\newcommand{\eenum}{\end{enumerate}}
\newcommand{\N}{\ensuremath{ \mathbf N }}

\newcommand{\Q}{\ensuremath{\mathbf Q}}
\newcommand{\R}{\ensuremath{\mathbf R}}

\newcommand{\C}{\ensuremath{\mathbf C}}

\newcommand{\mcc}{\ensuremath{ \mathcal C}}
\newcommand{\mcd}{\ensuremath{ \mathcal D}}

\newcommand{\mba}{\ensuremath{ \mathbf a}}
\newcommand{\mbb}{\ensuremath{ \mathbf b}}

\newcommand{\mbm}{\ensuremath{ \mathbf m}}
\newcommand{\mbo}{\ensuremath{ \mathbf 0}}

\newcommand{\mbw}{\ensuremath{ \mathbf w}}
\newcommand{\mbx}{\ensuremath{ \mathbf x}}
\newcommand{\mby}{\ensuremath{ \mathbf y}}
\newcommand{\mbz}{\ensuremath{ \mathbf z}}

\newcommand{\bmat}{\left(\begin{matrix}}
\newcommand{\emat}{\end{matrix}\right)}
\newcommand{\bsmallmat}{\left(\begin{smallmatrix}}
\newcommand{\esmallmat}{\end{smallmatrix}\right)}
\DeclareMathOperator{\qand}{\quad\text{and}\quad}
\DeclareMathOperator{\qqand}{\qquad\text{and}\qquad}

\title[Linear   and multiplicative polynomial equations]
{Linear equations and multiplicative polynomial equations in infinitely many variables}
\author{Melvyn B.   Nathanson}
\address{Lehman College (CUNY), Bronx, NY 10468}
\email{melvyn.nathanson@lehman.cuny.edu}

\author{David A. Ross}
\address{University of Hawaii at Manoa, Honolulu, HI 96822}
\email{ross@math.hawaii.edu}

\date{\today}



 \thanks{MBN supported in part by  PSC-CUNY Research Award Program grant 66197-00 54.}
\begin{document}

\begin{abstract}
This paper describes infinite sets of linear and polynomial equations 
in infinitely many variables with the property that the existence of a solution 
or even an approximate solution  
for every finite subset of the equations implies 
the existence of a solution for the infinite set of equations. 
\end{abstract}

\maketitle

\section{Finitely many implies infinite many}
In mathematics there are theorems asserting that, for certain classes of equations,  
if every finite subset of an infinite set of the equations has a  solution, 
then the infinite set of equations has a  solution.  
For example, if every finite subset of an infinite set of linear equations in $n$ variables 
(that is, equations of the form $\sum_{j=1}^n a_{i,j}x_j = b_i$) 
has a solution, 
then the  infinite set of linear equations in $n$ variables has a solution. 
For infinite sets of linear equations in infinitely many variables, in which each equation 
contains only finitely many variables 
(that is, equations of the form $\sum_{j=1}^{\infty} a_{i,j}x_j = b_i$ such that, for each $i$, 
we have  $a_{i,j} \neq 0$ for only finitely many $j$), 
a ``finitely solvable implies infinitely solvable'' theorem is also true.  
For a survey of related results, see Nathanson~\cite{nath24x}.

There are analogous results for linear equations that contain infinitely many variables.  
The pair  $(p,q)$  is \emph{H\" older conjugate}  
or, simply, \emph{conjugate} if  $p > 1$ and $q > 1$ are  real numbers such that 
\[
\frac{1}{p} + \frac{1}{q} = 1   
\]  
or if $p=1$ and $q = \infty$.  
For $1 \leq p < \infty$, let 
\[
\ell^p = \left\{ \mba = \left(a_j\right)_{j=1}^{\infty} : 
\|\mba\|_p = \left(\sum_{j=1}^{\infty} |a_j|^p \right)^{1/p}< \infty  \right\}
\]
and, for $p = \infty$, let 
\[
\ell^{\infty} = \left\{ \mba = \left(a_j\right)_{j=1}^{\infty} : 
\|\mba\|_{\infty} =\sup\{ |a_j| : j = 1,2,3,\ldots \} < \infty \right\}
\]
be the classical Lebesque spaces of
infinite sequences of real or complex numbers.

Recall H\" older's inequality:  If $(p,q)$ is a conjugate pair and 
\[
\mba = \left( a_j \right)_{j=1}^{\infty} \in \ell^p \qqand  
\mbx = \left(x_j \right)_{j=1}^{\infty} \in \ell^q 
\] 
then 
\[
\mba\mbx = \left( a_jx_j \right)_{j=1}^{\infty} \in \ell^1
\]
 and 
\[
\sum_{j=1}^{\infty} | a_jx_j  | = \left\|  \mba \mbx  \right\|_1 
\leq \left\|  \mba \right\|_p \left\|  \mbx  \right\|_q.
\]
Thus, the infinite series 
\[
 (\mba,\mbx) = \sum_{j=1}^{\infty} a_j x_j   
\] 
converges absolutely and 
it makes sense to ask if there is a solution $\mbx \in \ell^q$ of the linear equation 
in infinitely many variables 
\beq                             \label{MP:LinearEquation}
L(\mbx) =  (\mba,\mbx) = \sum_{j=1}^{\infty} a_jx_j = b.
\eeq

We retain the usual ambiguity between sequences of numbers and sequences of variables.  
An \emph{exact solution} or, simply, a \emph{solution}  of equation~\eqref{MP:LinearEquation} 
is a sequence 
$\mbx = (x_j)_{j=1}^{\infty} \in \ell^q$ 
such that $L(\mbx) = b$.  
Equation~\eqref{MP:LinearEquation} has an \emph{approximate solution} if, for every $\varepsilon > 0$, there exists $\mbx_{\varepsilon} \in \ell^q$ such that 
\[
\left| L(\mbx_{\varepsilon}) - b \right| < \varepsilon.
\]

Let $\N = \{1,2,3,\ldots\}$ be the set of positive integers 
and $\N^k$ the set of $k$-tuples of positive integers. 
Extending a classical result of F. Riesz~\cite[pp. 61--62]{ries13} 
 (see also Banach~\cite[p. 47]{bana87}),  
 Abian and Eslami~\cite{abia-esla82}  proved the following beautiful theorem.

\bt                \label{MP:theorem:Abian-Eslami}
Let $(p,q)$ be a conjugate pair with $p>1$ and $q>1$ 
and let $\mba_i =  \left( a_{i,j} \right)_{j=1}^{\infty} \in \ell^p$ 
for all $i \in \N$.  Consider the linear equation in infinitely many variables 
\[
L_i(\mbx) =  \sum_{j=1}^{\infty} a_{i,j} x_j = b_i. 
\]
Let $M > 0$. 
If, for every finite subset $S$ of \N, there is a sequence $\mbx_S \in \ell^q$ 
such that  $\|\mbx_S\|_q \leq M$ and $L_i(\mbx_S) = b_i$ for all $i \in S$, 
then there is a sequence $\mbx \in \ell^q$ 
such that  $\|\mbx\|_q \leq M$ and $L_i(\mbx) = b_i$ for all $i \in \N$. 
\et

Abian and Eslami proved their theorem only for $p > 1$ and for 
countably many equations 
with real coefficients and real solutions, but the result also holds for $(p,q) = (1,\infty)$ and 
for uncountably many equations with complex coefficients and complex solutions. 
In Appendix{~\ref{MP:appendix:Riesz} we prove Riesz's theorem and derive 
Theorem~\ref{MP:theorem:Abian-Eslami} from it. 

It is important to observe that Theorem~\ref{MP:theorem:Abian-Eslami} 
is not true without the condition 
that the norm of $\| \mbx_S \|_q$ is uniformly bounded by $M$ for all finite sets $S$.  
The following example is essentially due to Helly~\cite{hell21}.   

Let $\mba = (a_j)_{j=1}^{\infty} \in \ell^p$ with $a_j \neq 0$ for all $j \in \N$. 
For all $i \in\N$, define the sequence $\mba_i = (a_{i,j})_{j=1}^{\infty} \in \ell^p$ as follows:
\[
a_{i,j} = \begin{cases}
0 & \text{if $j < i$} \\
a_i & \text{if $j \geq i$.}
\end{cases}
\]
For all $i \in \N$, consider the  linear equation in infinitely many variables 
\[
L_i(\mbx) =  \sum_{j=1}^{\infty}  a_{i,j}x_j =  \sum_{j=i}^{\infty}  a_jx_j = 1.
\]
We obtain the ``triangular'' system of equations 
\[
\begin{array}{rl}
a_1x_1 + a_2x_2 + a_3 x_3 + \cdots + a_i x_i  + a_{i+1}x_{i+1} +\cdots  & = 1 \\
a_2x_2 + a_3 x_3 + \cdots + a_i x_i  + a_{i+1}x_{i+1} +\cdots  & = 1 \\
a_3 x_3 + \cdots + a_i x_i  + a_{i+1}x_{i+1} +\cdots  & = 1\\
\vdots & \\
a_i x_i  + a_{i+1}x_{i+1} +\cdots  & = 1\\
 a_{i+1}x_{i+1} +\cdots  & = 1\\
\end{array}
\]
For all $i \geq 1$, if  $L_i(\mbx) = L_{i+1}(\mbx) = 1$, then $a_ix_i =0$ and so $x_i = 0$ 
(because $a_i \neq 0$).  
It follows that the infinite set of equations has no solution.  
However, every finite subset of these equations is solvable.  

Here is one solution of the triangular system.  
For all $r \in \N$,  the sequence $\mbx_r = (x_{r,j})_{j=1}^{\infty} \in \ell^q$ defined by 
\[
x_{r,j} = \begin{cases}
1/a_r & \text{if $j = r$} \\
0& \text{if $j \neq r$}
\end{cases}
\]
is a solution of the finite set of equations $\{L_i(\mbx) = 1:i=1,\ldots, r\}$.  
We have 
\[
\|\mbx_r\|_q = \frac{1}{|a_r|} \qqand 
\lim_{r\rightarrow \infty} \|\mbx_r\|_q = \lim_{r\rightarrow \infty} \frac{1}{|a_r|} = \infty.  
\]
Thus, the sequence of solutions $(\mbx_r)_{r=1}^{\infty}$ is not uniformly bounded in $\ell^q$. 

More generally, if $\mbx_r$ is any solution of the first $r$ equations of the triangular system, then 
\[
1 = \sum_{j=r}^{\infty} a_jx_j = \sum_{j=1}^{\infty} a_{r,j}x_j  \leq \|  \mba_r\mbx_r\|_1 
\leq   \|  \mba_r \|_p \|\mbx_r\|_q
\]
and so 
\[
 \|\mbx_r\|_q \geq \frac{1}{  \|  \mba_r \|_p }.
\]
Because $\mba \in \ell^p$, we have  $\lim_{r\rightarrow \infty}   \|  \mba_r \|_p = 0$
and so  the sequence of solutions $(\mbx_r)_{r=1}^{\infty}$ is not uniformly bounded in $\ell^q$.

The following result extends Theorem~\ref{MP:theorem:Abian-Eslami} 
to approximately solvable systems of linear equations.  

\bt                            \label{MP:theorem:linear-approx}
Let $(p,q)$ be a conjugate pair.  
Let  $I$ be an infinite set.   For all $i \in I$, let 
$\mba_i  = \left( a_{i,j} \right)_{j=1}^{\infty} \in \ell^p$      
and consider the linear equation in infinitely many variables 
\[
L_i(\mbx) =  \sum_{j=1}^{\infty} a_{i,j} x_j =  b_i.  
\]
Let $M > 0$. 
If, for every finite subset $S$ of I and every $\varepsilon > 0$, 
the finite set of linear inequalities 
\[
\left\{ |L_i(\mbx) -b_i | \leq \varepsilon:i \in S \right\} 
\]
has a solution $\mbx_{S,\varepsilon} \in \ell^q$ 
with $\| \mbx_{S,\varepsilon} \|_q \leq M$, 
then the infinite set of linear equations 
\[
\left\{ L_i(\mbx) = b_i:i \in I \right\} 
\] 
has an exact solution $\mbx \in \ell^q$ with $\| \mbx \|_q \leq M$. 
\et

In this paper we generalize Theorem~\ref{MP:theorem:linear-approx} to certain infinite 
sets of polynomial equations in infinitely many variables and prove 
 that the existence of norm-bounded or even sequentially bounded
approximate solutions to all finite subsets of the set of polynomial equations is sufficient 
to guarantee the existence of an exact solution to the infinite set (Theorems~\ref{MP:theorem:multiplicative-approx} and~\ref{MP:theorem:multiplicative-approx-m}).  
The theorem on linear equations is a  special case of the polynomial results.

Our theorems apply to polynomials with coefficients and solutions 
in the field of real numbers and also in the field of complex numbers.  
However, the ``finitely many implies infinitely many'' paradigm is not true 
for all subfields of the complex numbers.
For example, the paradigm fails for the field \Q\ of rational numbers. 
Let $(b_i)_{i\in I}$ be an infinite sequence of irrational numbers and consider 
the infinite set of linear equations $\{ x_i = b_i:i \in I\}$. 
For every $\varepsilon > 0$, every finite subset of the equations and, indeed, the infinite set  
of  equations have approximate solutions in \Q, but neither the infinite set 
nor any nonempty finite subset of the equations has an exact solution in \Q.

\section{Polynomials in infinitely many variables} 
Let $D$ be a positive integer, let $d \in \{1,2,\ldots, D\}$, and let 
\[
 D < q \leq \infty.
\] 
If $q < \infty$, then 
\[
\frac{1}{q/(q-d)} + \frac{1}{q/d} = 1 
\]
and so the real numbers $q/(q-d) > 1$ and $q/d > 1$ form a conjugate pair. 
For $q = \infty$, we define 
\[
\frac{q}{(q-d)} = 1 \qqand  \frac{q}{d} = \infty.  
\]
If 
\[
\mba= \left( a_j \right)_{j=1}^{\infty} \in \ell^{q/(q-d)} 
\]
and 
\[
\mbx = \left(x_j \right)_{j=1}^{\infty} \in \ell^q 
\] 
then 
\[
\mbx^d = \left(x_j^d \right)_{j=1}^{\infty} \in \ell^{q/d}.  
\]
From H\" older's inequality we obtain  
\[
\mba\mbx^d  = \left(a_j x_j^d \right)_{j=1}^{\infty} \in \ell^1 
\]
 and 
\[
 \sum_{j=1}^{\infty} \left| a_j x_j^d  \right| = \left\| \mba\mbx^d \right\|_1 
  \leq  \left\| \mba \right\|_{q/(q-d)}  \left\|\mbx^d \right\|_{q/d} < \infty. 
\]
Thus, the infinite series 
\[
 \left( \mba,\mbx^d \right) = \sum_{j=1}^{\infty} a_j x_j^d  
\] 
converges absolutely.

For all $k \in \{1,2,\ldots, D\}$,  let $\mcd_k$  be the set of all $k$-tuples of positive integers 
whose sum is at most $D$, that is,  
\beq                                                  \label{MP:Dk}
\mcd_k =  \left\{  (d_1,\ldots, d_k) \in \N^k:  d_1 + d_2 + \cdots + d_k \leq D \right\}.
\eeq
This is a finite set of cardinality  $\binom{D}{k}$.  
The set $\bigcup_{k=1}^D \mcd_k$ has cardinality $\sum_{k=1}^D \binom{D}{k} = 2^D -1$. 
 For all  
 \[
 \Delta = (d_1, d_2, \ldots, d_k) \in \mcd_k
 \]
  and 
 \[
 J = (j_1, j_2, \ldots, j_k) \in \N^k   
\]
we define the monomial 
\[
x_J^{\Delta} = x_{j_1}^{d_1}  x_{j_2}^{d_2}  \cdots  x_{j_k}^{d_k}    
\] 
of degree 
\[
|\Delta| = d_1 + d_2 + \cdots + d_k \leq D.
\]

For example, let  $D=3$ and $k \in \{1,2,3\}$.    
There are $7 = 2^3-1$ sets of monomials.
For $k=1$ and $J = (j_1) \in \N^1$, we have 
\[ 
\Delta  \in \mcd_1  = \{ (1), (2), (3) \}  \qand 
x_J^{\Delta} \in \left\{ x_{j_1} , \ x_{j_1}^2 , \ x_{j_1}^3 \right\}.
\]
For $k=2$ and $ J = (j_1,j_2) \in \N^2$, we have 
\[
\Delta  \in \mcd_2   = \{ (1,1), (1,2), (2,1) \}
 \qand x_J^{\Delta} \in \left\{ x_{j_1} x_{j_2}, \ x_{j_1}x_{j_2}^2 , \ x_{j_1}^2x_{j_2} \right\}.
\] 
For $k=3$ and $ J = (j_1,j_2,j_3) \in \N^3$, we have 
\[
\Delta  \in \mcd_3  = \{ (1,1,1) \} \qand x_J^{\Delta} = x_{j_1} x_{j_2} x_{j_3}.
\]

Monomials do not have a unique representation in the form $x_J^{\Delta}$. 
For example, with $D=3$, let  $k=2$ and $J = (1,1) \in \N^2$.  
If $\Delta = (1,2) \in \mcd_2$, then $x_J^{\Delta} = x_1x_1^2  = x_1^3$. 
If $\Delta = (2,1)\in \mcd_2$, then $x_J^{\Delta} = x_1^2x_1  = x_1^3$. 
If $D=3$, $k=1$, $J= (1) \in \N^1$, and $\Delta = (3) \in \mcd_1$, then 
$x_J^{\Delta} = x_1^3$. 
Note that we do not need unique representation of monomials, 
but we would have unique representation if we considered only $k$-tuples 
$(j_1, j_2, \ldots, j_k) \in \N^k$ such that $j_1 < j_2 < \cdots < j_k$.)  

We consider only polynomials with zero constant term.  In this paper, ``polynomial'' 
means polynomial with zero constant term.
By moving constants from one side of an equation to the other side, 
it is always possible to transform an equation that includes polynomials with nonzero constant terms 
to an equivalent equation in which no polynomial has a nonzero constant term.

A \emph{polynomial of degree at most $D$ in infinitely many variables} $\mbx = (x_j)_{j=1}^{\infty}$ 
is a formal power series of the form 
\beq                     \label{MP:polynomial}
P(\mbx) = \sum_{k=1}^D \sum_{ \Delta \in \mcd_k} \sum_{J \in \N^k} 
 a_{\Delta,J} \  x_J^{\Delta} 
\eeq
with real or complex coefficients $a_{\Delta,J}$.  
For example, a polynomial of degree at most 3 in infinitely many variables is of the form 
\begin{align*} 
P(\mbx) 
& =  \sum_{j_1=1}^{\infty} a_{(1),(j_1)} x_{j_1} +  \sum_{j_1=1}^{\infty} a_{(2),(j_1)} x_{j_1}^2
 +  \sum_{j_1=1}^{\infty} a_{(3),(j_1)} x_{j_1}^3 \\ 
 & +   \sum_{j_1=1}^{\infty} \sum_{j_2 =1}^{\infty} a_{(1,1), (j_1,j_2)} x_{j_1}  x_{j_2} 
 + \sum_{j_1=1}^{\infty} \sum_{j_2 =1}^{\infty} a_{(1,2), (j_1,j_2)} x_{j_1}  x_{j_2}^2 \\
& \quad  +  \sum_{j_1=1}^{\infty} \sum_{j_2 =1}^{\infty} a_{(2,1), (j_1,j_2)} x_{j_1}^2  x_{j_2} 
+   \sum_{j_1=1}^{\infty} \sum_{j_2 =1}^{\infty} \sum_{j_3 =1}^{\infty}  a_{(1,1,1), (j_1,j_2,j_3)} x_{j_1}  x_{j_2} x_{j_3}.
\end{align*}

\section{Multiplicative polynomials} 
Let $D$ be a positive integer and let $1 \leq D < q \leq \infty$. 
Let
\[
P(\mbx) = \sum_{k=1}^D \sum_{ \Delta \in \mcd_k} \sum_{J \in \N^k} 
 a_{\Delta,J} \  x_J^{\Delta} 
\]
be a polynomial of degree at most $D$ in infinitely many variables $\mbx = (x_j)_{j=1}^{\infty}$  
with coefficients $a_{\Delta,J}$.
We consider polynomials whose coefficients $a_{\Delta,J}$ are \emph{multiplicative} 
in the following sense: For all $d \in \{1,2,\ldots, D\}$  there is a sequence  
\[
\mba_d  = \left( a_{d,j} \right)_{j=1}^{\infty} \in \ell^{q/(q-d)}
\]
such that,
if 
\[
\Delta = (d_1, d_2, \ldots, d_k) \in \mcd_k  \qqand   J = (j_1, j_2, \ldots, j_k) \in \N^k 
\]
then 
\[
a_{\Delta,J} = a_{(d_1, d_2, \ldots, d_k),(j_1, j_2, \ldots, j_k)}  = a_{d_1,j_1} a_{d_2,j_2} \cdots a_{d_k,j_k}. 
\]
This gives the monomial factorization
\begin{align*}
a_{\Delta,J}  x_J^{\Delta}  
& = \left( a_{d_1,j_1} a_{d_2,j_2} \cdots a_{d_k,j_k} \right) 
\left(  x_{j_1}^{d_1}  x_{j_2}^{d_2}  \cdots  x_{j_k}^{d_k} \right) 
 \\
& = \left( a_{d_1,j_1} x_{j_1}^{d_1} \right) 
\left( a_{d_2,j_2} x_{j_2}^{d_2} \right) 
\cdots \left( a_{d_k,j_k}   x_{j_k}^{d_k} \right)  
\end{align*}
and so 
\begin{align*} 
 \sum_{J \in \N^k}  a_{\Delta,J} \  x_J^{\Delta} 
 & =  \sum_{ (j_1,j_2,\ldots, j_j) \in \N^k} 
  \left( a_{d_1,j_1} x_{j_1}^{d_1} \right) 
\left( a_{d_2,j_2} x_{j_2}^{d_2} \right) 
\cdots \left( a_{d_k,j_k}   x_{j_k}^{d_k} \right)\\
& =  \left( \sum_{j_1=1}^{\infty} a_{d_1,j_1} x_{j_1}^{d_1} \right) 
\left(\sum_{j_2=1}^{\infty} a_{d_2,j_2} x_{j_2}^{d_2} \right) \cdots 
\left( \sum_{j_k=1}^{\infty} a_{d_k,j_k}   x_{j_k}^{d_k}  \right)    \\
& =     
\left( \mba_{d_1}, \mbx^{d_1 } \right) \left( \mba_{d_2}, \mbx^{d_2 } \right) 
\cdots \left( \mba_{d_k}, \mbx^{d_k } \right). 
\end{align*}
For $\mbx = (x_j)_{j=1}^{\infty} \in \ell^q$, the rearrangement is justified 
by the absolute convergence of the $k$ infinite series 
$\left( \mba_{d_1}, \mbx^{d_1 } \right), \left( \mba_{d_2}, \mbx^{d_2 } \right),
\ldots, \left( \mba_{d_k}, \mbx^{d_k } \right)$.  
We obtain  
\begin{align*}
P(\mbx) 
& = \sum_{k=1}^D \sum_{ \Delta \in \mcd_k} \sum_{J \in \N^k}  a_{\Delta,J} \  x_J^{\Delta} \\
& =   \sum_{k=1}^D \sum_{ (d_1, d_2, \ldots, d_k)  \in \mcd_k}  
\left( \sum_{j_1=1}^{\infty} a_{d_1,j_1} x_{j_1}^{d_1} \right) 
\left(\sum_{j_2=1}^{\infty} a_{d_2,j_2} x_{j_2}^{d_2} \right) \cdots 
\left( \sum_{j_k=1}^{\infty} a_{d_k,j_k}   x_{j_k}^{d_k}  \right)    \\
& =      \sum_{k=1}^D \sum_{ (d_1, d_2, \ldots, d_k)  \in \mcd_k}     
\left( \mba_{d_1}, \mbx^{d_1 } \right)\left( \mba_{d_2}, \mbx^{d_2 } \right) 
\cdots \left( \mba_{d_k}, \mbx^{d_k } \right). 
\end{align*}
This is a finite sum, and so the series $P(\mbx)$ converges absolutely 
for all $\mbx \in \ell^q$.

Let $D$ be a positive integer and let $1 \leq D < q \leq \infty$.  
A \emph{$q$-multiplicative polynomial of degree at most $D$}\index{multiplicative polynomial} 
or, simply, a \emph{multiplicative polynomial}, 
is a polynomial of degree at most $D$ in infinitely many variables of the form 
\[
P(\mbx) =    \sum_{k=1}^D \sum_{ (d_1, d_2, \ldots, d_k)  \in \mcd_k}  
\left( \mba_{d_1}, \mbx^{d_1 } \right)\left( \mba_{d_2}, \mbx^{d_2 } \right) 
\cdots \left( \mba_{d_k}, \mbx^{d_k } \right)
\]
where 
\[
\mba_{d}  = \left( a_{d,j} \right)_{j=1}^{\infty} \in \ell^{q/(q-d)}
\]
for all $d \in \{1,\ldots, D\}$.  The infinite series converges for all $\mbx \in \ell^q$. 

A \emph{multiplicative polynomial equation} is an equation of the form   
\[
P(\mbx) = b
\]
where $P(\mbx)$ is a multiplicative polynomial.  
We ask if this equation has a solution $\mbx \in \ell^q$.

The analogue of Theorem~\ref{MP:theorem:Abian-Eslami} 
for linear equations is the following  ``finitely many implies infinitely many'' solvability result 
for multiplicative polynomial equations.

\bt                                            \label{MP:theorem:multiplicative}
Let $D$ be a positive integer and let $1 \leq D < q \leq \infty$. 
Let $I$ be an infinite set.  
For all $i \in I$ and $d \in \{1,2,\ldots, D\}$, let 
\[
\mba_{i,d}  = \left( a_{i,d,j} \right)_{j=1}^{\infty} \in \ell^{q/(q-d)}.
\]
For $\Delta = (d_1,d_2,\ldots, d_k) \in \mcd_k$ and $J = (j_1, j_2, \ldots, j_k) \in \N^k$, 
let 
 \[
a_{i,\Delta,J} = a_{i,d_1,j_1} a_{i,d_2,j_2} \cdots a_{i,d_k,j_k}.
\]
For all $i \in I$, consider the multiplicative polynomial equation 
\begin{align*}
P_i(\mbx) 
& = \sum_{k=1}^D \sum_{ \Delta \in \mcd_k} \sum_{J \in \N^k}  a_{i,\Delta,J} \  x_J^{\Delta} \\
& =     \sum_{k=1}^D \sum_{ (d_1, d_2, \ldots, d_k)  \in \mcd_k}     
\left( \mba_{i,d_1}, \mbx^{d_1 } \right)\left( \mba_{i,d_2}, \mbx^{d_2 } \right) 
\cdots \left( \mba_{i,d_k}, \mbx^{d_k } \right) \\ 
& = b_i. 
\end{align*}
 Let $M > 0$. 
If, for every finite subset $S$ of I, the finite set of polynomial equations 
$\{ P_i(\mbx) = b_i:i \in S \}$ 
has a solution $\mbx_S \in \ell^q$ with $\| \mbx_S \|_q \leq M$, 
then the infinite set of polynomial equations $\{P_i(\mbx) = b_i:i \in I\}$ 
has a solution $\mbx \in \ell^q$ with $\| \mbx \|_q \leq M$. 
\et

\section{Approximate finite implies exact infinite}

The set $\{P_i(\mbx) = b_i:i \in S \}$ of polynomial equations in infinitely many 
variables $\mbx = (x_j)_{j=1}^{\infty}$ has an  \emph{approximate solution} if, 
for every $\varepsilon > 0$, there exists a sequence $\mbx_{\varepsilon}$ such that 
$|P_i(\mbx_{\varepsilon}) - b_i| \leq \varepsilon$ for all $i \in S$.
Theorem~\ref{MP:theorem:multiplicative-approx} states that  an infinite set 
of multiplicative polynomial equations  has an exact solution 
if every finite subset of the equations has a 
norm-bounded approximate solution.   
This is the main result of this paper and immediately implies 
Theorems~\ref{MP:theorem:Abian-Eslami}, \ref{MP:theorem:linear-approx}, 
and~\ref{MP:theorem:multiplicative}.

The theorem and the lemmas in this section 
are valid in both the real and complex fields.  

\bt                            \label{MP:theorem:multiplicative-approx} 
Let $D$ be a positive integer and let 
$1 \leq D < q \leq \infty$.   For all $k \in \{1,2,\ldots, D\}$, let 
\[
\mcd_k =  \left\{  (d_1,\ldots, d_k) \in \N^k:  d_1 + d_2 + \cdots + d_k \leq D \right\}.
\]
Let $I$ be an infinite set and let $\mbx = (x_j)_{j=1}^{\infty}$.  
For all $i \in I$ and $d \in \{1,2,\ldots, D\}$, let 
\[
\mba_{i,d}  = \left( a_{i,d,j} \right)_{j=1}^{\infty} \in \ell^{q/(q-d)}   
\] 
and  
\[
\left(\mba_{i,d},\mbx^d \right) = \sum_{j=1}^{\infty} a_{i,d,j}x_j^d. \]
For $\Delta = (d_1,d_2,\ldots, d_k) \in \mcd_k$ and $J = (j_1, j_2, \ldots, j_k) \in \N^k$, 
define the multiplicative coefficient 
 \[
a_{i,\Delta,J} = a_{i,d_1,j_1} a_{i,d_2,j_2} \cdots a_{i,d_k,j_k}.
\]  
For all $i \in I$, the finite set of sequences $\{\mba_{i,d} \}_{d=1}^D$ 
determines the multiplicative polynomial 
\begin{align*}
P_i(\mbx) 
& = \sum_{k=1}^D \sum_{ \Delta \in \mcd_k} \sum_{J \in \N^k}  a_{i,\Delta,J} \  x_J^{\Delta} \\
& =     \sum_{k=1}^D \sum_{ (d_1, d_2, \ldots, d_k)  \in \mcd_k}     
\left( \mba_{i,d_1}, \mbx^{d_1 } \right)\left( \mba_{i,d_2}, \mbx^{d_2 } \right) 
\cdots \left( \mba_{i,d_k}, \mbx^{d_k } \right). 
\end{align*} 
 Let $M > 0$ and $\mbb = (b_i)_{i\in I}$. 
If, for every $\varepsilon > 0$ and every finite subset $S$ of $I$, 
the finite set of polynomial inequalities 
\[
\left\{ |P_i(\mbx) -b_i | \leq \varepsilon:i \in S \right\} 
\]
has a solution $\mbx_{S,\varepsilon} \in \ell^q$ with $\| \mbx_{S,\varepsilon} \|_q \leq M$, 
then the infinite set of polynomial equations 
\[
\left\{ P_i(\mbx) = b_i:i \in I \right\} 
\] 
has an exact solution $\mbx \in \ell^q$ with $\| \mbx \|_q \leq M$. 
\et

We begin with two results whose statements and 
proofs are valid in both the real and complex cases.  
 
For $M > 0$, the closed interval  $[-M,M] = \{x\in \R: |x| \leq M\}$ 
and the closed ball $B_M = \{x\in \C: |x|  \leq M\}$ are compact.  
For polynomial equations in \R\ we use the compact topological space  
\[
\Omega = \prod_{j=1}^{\infty} [-M,M]
\]
and for polynomial equations in \C\ we use the compact  topological  space  
\[
\Omega = \prod_{j=1}^{\infty} B_M.
\]
Let 
\[
X_{q,M} = \left\{ \mbx \in \ell^q: \| \mbx \|_q \leq M\right\}.  
\]

\bl               \label{MP:lemma:compact}
The set $X_{q,M}$ is a compact subset of the topological space $\Omega$. 
\el

\begin{proof} 
If $\mbx = (x_j)_{j=1}^{\infty} \in \ell^q$ and $\|\mbx\|_q \leq M$, then $|x_j| \leq M$ for all
$j \in \N$ and so $X_{q,M} \subseteq \Omega$.
Because $\Omega$ is compact, it suffices to prove that $X_{q,M}$ is closed 
in the product topology on $\Omega$, or equivalently, that the complement of $X_{q,M}$ 
is open.  
The complement of $X_{q,M}$ is the set 
\begin{align*}
Y_{q,M} & = \Omega \setminus X_{q,M} \\ 
& = \left\{ \mby \in \Omega: \mby \notin \ell^q  \right\} 
\bigcup  \left\{ \mby \in \Omega: \mby \in \ell^q \text{ and } \|\mby\|_q > M \right\} \\ 
& =  \left\{ \mby = (y_j)_{j=1}^{\infty} : \sum_{j=1}^{\infty} |y_j|^q > M^q \right\}.  
\end{align*}

Let $\mby   = (y_j)_{j=1}^{\infty} \in Y_{q,M}$.  
There exist $\varepsilon > 0$ and $N \in \N$ such that  
\[
\sum_{j=1}^{N} |y_j|^q > M^q + \varepsilon.
\]
Because the function $f(t) = |t|^q$ is  continuous, there exists $\delta > 0$ such 
that  $|t-y_j| < \delta$ implies   
\[
 |t|^q > |y_j|^q -  \frac{\varepsilon}{2N} 
\]
for all $j \in \{1,2,\ldots, N\}$. 
The set 
\[
U = \left\{ \mbz  = (z_j)_{j=1}^{\infty} \in \Omega: |z_j - y_j| < \delta 
\text{ for all } j \in \{1,2,\ldots, N\}  \right\}
\]
is an open neighborhood of $\mby$ in $\Omega$.  
For all $\mbz  = (z_j)_{j=1}^{\infty} \in U$ we have  
\begin{align*}
\sum_{j=1}^{\infty} |z_j|^q & \geq 
\sum_{j=1}^{N} |z_j|^q \ >  \sum_{j=1}^{N} \left(  |y_j|^q -  \frac{\varepsilon}{2N}  \right) \\
& =  \sum_{j=1}^{N}  |y_j|^q - \frac{\varepsilon}{2}   > M^q + \frac{\varepsilon}{2} \\
&  > M^q
\end{align*}
and so  $\mbz \in Y_{q,M}$.  
It follows that $U \subseteq Y_{q,M}$ and so 
$Y_{q,M}$ is an open subset of $\Omega$.
This completes the proof.  
\end{proof}

\bl                         \label{MP:lemma:continuous-product} 
Let  $d$ be a positive integer, $d < q \leq \infty$, and let $p = q/(q-d)$.  
Let $\mba = (a_j)_{j=1}^{\infty} \in \ell^p$.
Let $M > 0$.  
The function $f_d$ on $X_{q,M}$ defined by 
\[
f_d(\mbx) = (\mba,\mbx^d) = \sum_{j=1}^{\infty} a_j x_j^d 
\] 
is continuous with respect to the product topology on $X_{q,M} $ 
as a subspace of $\Omega$. 
\el

\begin{proof}
If $\mbx \in X_{q,M}$, then $\mbx^d \in \ell^{q/d}$.  
Because $(p,q/d)$ is a conjugate pair, 
the infinite series $\sum_{j=1}^{\infty} a_j x_j^d $ converges absolutely and the 
function $f_d(\mbx)$ is well-defined.   

Let $U$ be an open subset of \R\ or \C.  
We shall prove $f_d^{-1}(U)$ is open in $X_{q,M}$ 
as a subspace of $\Omega$.

Let $\mbx = (x_j)_{j=1}^{\infty} \in f_d^{-1}(U)$.  
Because $f_d(\mbx)  \in U$, there exists $\varepsilon > 0$ such that $U$ 
contains the open set
\[
\left\{ t : |t - f_d(\mbx) | < \varepsilon \right\}.
\]  
Because $\mba \in \ell^p$, the series 
$\sum_{j=1}^{\infty} |a_j|^p$ converges and there is an integer $N_{\mba}$ such that 
\[
\sum_{j=N+1}^{\infty} |a_j|^p < \left( \frac{\varepsilon}{3M} \right)^p 
\]
for all $N \geq N_{\mba}$.
Because $\mbx \in \ell^q$, the series 
$\sum_{j=1}^{\infty} \left| x_j \right|^q$ converges 
and there is an integer $N_{\mbx}$ such that 
\[
\sum_{j=N+1}^{\infty} |x_j|^q < \left( \frac{\varepsilon}{3M} \right)^q 
\]
for all $N \geq N_{\mbx}$.
Choose $N \geq \max(N_{\mba},N_{\mbx})$ and let $\delta > 0$ satisfy 
\[
\delta \sum_{j=1}^N |a_j| < \frac{\varepsilon}{3}. 
\]
Choose $\delta' > 0$ such that $ |y_j - x_j| < \delta'$ implies 
$ |y_j^d - x_j^d| < \delta$ for all  $j \in \{1,2,\ldots,N\}$. 
The set 
\[
V' = \left\{ \mby = (y_j)_{j=1}^{\infty} \in \Omega : |y_j - x_j| < \delta'
\text{ for all } j \in \{1,2,\ldots,N\} \right\}
\]
is an  open neighborhood of \mbx\ in the topological space $\Omega$ and so 
\begin{align*}
V & = V' \cap X_{q,M} \\ 
& = \left\{ \mby = (y_j)_{j=1}^{\infty} \in X_{q,M} : |y_j - x_j| < \delta'  
\text{ for all } j \in \{1,2,\ldots,N\} \right\}
\end{align*}
is an open neighborhood of \mbx\ in $X_{q,M}$.  
We shall prove that $V \subseteq f_d^{-1}(U)$. 

Let $\mby = (y_j)_{j=1}^{\infty} \in V$.  We define 
\[
\tilde{x_j} = \begin{cases}
0 & \text{if $j \leq N$} \\
x_j & \text{if $j \geq N+1$} 
\end{cases}
\]
and 
\[
\tilde{y_j} = \begin{cases}
0 & \text{if $j \leq N$} \\
y_j & \text{if $j \geq N+1$.} 
\end{cases}
\]
Then 
\[
 \tilde{\mbx} = ( \tilde{x_j} )_{j=1}^{\infty} \in \ell^q \qqand 
\tilde{\mby} = ( \tilde{y_j} )_{j=1}^{\infty} \in \ell^q  
\] 
and  
\[
\| \tilde{\mbx} \|_q \leq \| \mbx  \|_q \leq M
\qqand 
\| \tilde{\mby} \|_q \leq \| \mby  \|_q \leq M.
\]
Thus, $\tilde{\mbx} \in X_{q,M}$ and $\tilde{\mby} \in X_{q,M}$.
We define 
\[
\tilde{a_j} = \begin{cases}
0 & \text{if $j \leq N$} \\
a_j & \text{if $j \geq N+1$.} 
\end{cases}
\]
Then  
\[
 \tilde{\mba} = ( \tilde{a_j} )_{j=1}^{\infty} \in \ell^p 
 \]
and  
 \[
\| \tilde{\mba} \|_p \leq \frac{\varepsilon}{3M}.
 \]
We have 
\begin{align*}
 f_d(\mby) - f_d(\mbx) 
 & = (\mba, \mby^d) -  (\mba, \mbx^d ) \\
& = \sum_{j=1}^{\infty} a_jy_j^d  -  \sum_{j=1}^{\infty} a_j x_j^d  
 = \sum_{j=1}^{\infty} a_j \left( y_j^d  - x_j^d \right) \\
& =   \sum_{j=1}^{N} a_j(y_j^d - x_j^d)  +   \sum_{j=N+1}^{\infty} a_j(y_j^d - x_j^d). 
\end{align*}
Applying the triangle, H\"older, and Minkowski inequalities, we obtain  
\begin{align*}
\left| f_d(\mby) - f_d(\mbx) \right| 
& \leq  \sum_{j=1}^{N} |a_j| \ |y_j^d - x_j^d|  +   \sum_{j=N+1}^{\infty} |a_j(y_j^d - x_j^d) | \\ 
& \leq \delta \sum_{j=1}^{N} |a_j|  + \left\| \left( \tilde{\mba}, \tilde{\mby}  - \tilde{\mbx}  \right) \right\|_1  \\ 
& < \frac{\varepsilon}{3}  
+ \left\| \tilde{\mba} \right\|_p  \left\| \tilde{\mby}  - \tilde{\mbx} \right\|_q  \\
& \leq \frac{\varepsilon}{3} +  \left\| \tilde{\mba} \right\|_p 
\left( \left\|   \tilde{\mby}   \right\|_q +  \left\|  \tilde{\mbx} )  \right\|_q \right) \\
& < \frac{\varepsilon}{3} + \frac{\varepsilon}{3M} (M+M) \\
& = \varepsilon. 
\end{align*}
Therefore, $f_d(\mby) \in U$ and  $V \subseteq f_d^{-1}(U)$.  
Thus, the set $ f_d^{-1}(U)$ is open in $X_{q,M}$ and the function $f_d(\mbx) = (\mba, \mbx^d)$ 
is continuous on $X_{q,M}$. 
This completes the proof.  
\end{proof}

We can now prove Theorem~\ref{MP:theorem:multiplicative-approx}.

\begin{proof}
For all $i \in I$ we have the multiplicative polynomial 
\[
P_i(\mbx)  =     \sum_{k=1}^D \sum_{ (d_1, d_2, \ldots, d_k)  \in \mcd_k}     
\left( \mba_{i,d_1}, \mbx^{d_1 } \right)\left( \mba_{i,d_2}, \mbx^{d_2 } \right) 
\cdots \left( \mba_{i,d_k}, \mbx^{d_k } \right).  
\] 
By Lemma~\ref{MP:lemma:compact}, the set $X_{q,M}$ is a compact subset of $\Omega$.   
By Lemma~\ref{MP:lemma:continuous-product}, the 
functions $\left( \mba_{i,d_i}, \mbx^{d_i } \right)$ 
are continuous on $X_{q,M} $ for all $d_i \in \{1,\ldots, D\}$.  
Finite sums of finite products of continuous functions are continuous, 
and so the multiplicative polynomials 
$P_i(\mbx)$ are continuous functions on $X_{q,M}$  for all $i \in I$. 
It follows that, for all $i \in I$ and $\varepsilon > 0$, the approximation set  
\[
F_{i,\varepsilon} = \left\{ \mbx \in X_{q,M}: | P_i(\mbx) - b_i | \leq \varepsilon  \right\}
\]
is a  closed subset of the compact set $X_{q,M} $.  For every finite subset $S$ of $I$,
the set of  polynomial inequalities 
\[
\left\{ |P_i(\mbx) -b_i | \leq \varepsilon:i \in S \right\} 
\]
has a solution $\mbx_{S,\varepsilon} \in \ell^q$ with $\| \mbx_{S,\varepsilon} \|_q \leq M$ 
and so the set of closed sets $\{F_{i,\varepsilon} : i \in I \text{ and } \varepsilon > 0 \}$ 
has the finite intersection property.  Therefore,  
\[
\bigcap_{\substack{i \in I \\ \varepsilon > 0}} F_{i,\varepsilon}  
= \left\{  \mbx \in X_{q,M}: P_i(\mbx) = b_i  
\text{ for all $i \in I$} \right\}
\] 
is nonempty and the infinite set of polynomial equations has an exact solution.    
This completes the proof Theorem~\ref{MP:theorem:multiplicative-approx}. 
\end{proof} 

Theorem~\ref{MP:theorem:linear-approx} for infinitely many linear equations 
in infinitely many variables is the special case of 
Theorem~\ref{MP:theorem:multiplicative-approx} with $D=1$.  
Theorem~\ref{MP:theorem:multiplicative-approx} implies Theorem~\ref{MP:theorem:multiplicative}  
and Theorem~\ref{MP:theorem:linear-approx} implies Theorem~\ref{MP:theorem:Abian-Eslami}.  

We have the following refinement of Theorem~\ref{MP:theorem:multiplicative-approx}.

\bt                            \label{MP:theorem:multiplicative-approx-m} 
With the hypotheses of Theorem~\ref{MP:theorem:multiplicative-approx}, let 
 \[
 \mbm = (m_j)_{j=1}^{\infty} \in \ell^q \qquad \text{ with $m_j \geq 0$ for all $j \in \N$. }
 \]
If, for every $\varepsilon > 0$ and every finite subset $S$ of I, 
the finite set of polynomial inequalities 
\[
\left\{ |P_i(\mbx) -b_i | \leq \varepsilon:i \in S \right\} 
\]
has a solution 
$\mbx_{S,\varepsilon} = \left( x_{S,\varepsilon,j}  \right)_{j=1}^{\infty}\in \ell^q$ 
with $| x_{S,\varepsilon,j} | \leq m_j$ for all $j \in \N$, 
then the infinite set of polynomial equations 
\[
\left\{ P_i(\mbx) = b_i:i \in I \right\} 
\] 
has an exact solution $\mbx = \left( x_j \right)_{j=1}^{\infty}\in \ell^q$ 
with $| x_j | \leq m_j$ for all $j \in \N$. 
\et

\begin{proof}
Let $M = \|\mbm\|_q$.  For polynomial equations in \R\ we use the compact space  
\[
X_{\mbm} = \prod_{j=1}^{\infty} [-m_j,m_j]
\]
and for polynomial equations in \C\ we use the compact space  
\[
X_{\mbm} = \prod_{j=1}^{\infty} B_{m_j}.
\]
If $\mbx \in X_{\mbm}$, then $\mbx \in \ell^q$ and $\|\mbx\|_q \leq M$, and so 
\[
X_{\mbm} \subseteq X_{q,M} =  \left\{ \mbx \in \ell^q: \|\mbx\|_q \leq M\right\}.
\]
By Lemma~\ref{MP:lemma:compact}, the set $X_{q,M}$ is a compact subset of $\Omega$. 
By Lemma~\ref{MP:lemma:continuous-product}, 
for all $d \in \{1,\ldots, D\}$ and for all $\mba \in \ell^{q/(q-d)}$, 
the  function $f_d(\mbx) = (\mba,\mbx^d)$ 
is continuous on $X_{q,M}$ and so its restriction to 
the compact set $X_{\mbm}$ is continuous.  
It follows that the multiplicative polynomials $P_i(\mbx)$ are 
continuous functions on $X_{\mbm}$ for all $i \in I$, 
and so  the approximation set  
\[
F_{i,\varepsilon} = \left\{ \mbx \in X_{\mbm}: | P_i(\mbx) - b_i | \leq \varepsilon  \right\}
\]
is a nonempty closed subset of the compact set $X_{\mbm} $.  
These sets have the finite intersection property and so  
\[
\bigcap_{\substack{i \in I \\ \varepsilon > 0}} F_{i,\varepsilon}  
= \left\{  \mbx \in X_{\mbm}: P_i(\mbx) = b_i  
\text{ for all $i \in I$} \right\} 
\] 
is nonempty.  
This completes the proof 
\end{proof}

\section{An application to Dirichlet series}

Here is an example of a multiplicative polynomial.   Let $D$ be a positive integer and let 
 $D < q \leq \infty$.  For  $d \in \{1,\ldots, D\}$, 
 let ${s_d} \in \C$ with $\Re({s_d}) > (q-d)/q$.
 Then 
 \[
\Re\left( \frac{qs_d}{q-d} \right)  > 1 
 \]
 and so the Riemann zeta function 
\[
\zeta\left( \frac{q{s_d}}{q-d}\right)  = \sum_{j=1}^{\infty} \frac{1}{j^{q{s_d}/(q-d)} } 
= \sum_{j=1}^{\infty} \left( \frac{1}{j^{s_d}}\right)^{q/(q-d)} 
\]
converges absolutely.  Consider the sequence 
\[
\mba_d  = \left( a_{d,j} \right)_{j=1}^{\infty} = \left( \frac{1}{j^{s_d}} \right)_{j=1}^{\infty} \in \ell^{q/(q-d)}    
\]
for all $d \in \{1,\ldots, D\}$, where 
\[
a_{d,j} =  \frac{1}{j^{s_d}} 
\]
and the polynomial 
\[
(\mba_d, \mbx^d)   = \sum_{j=1}^{\infty} a_{d,j} x_j^d = \sum_{j=1}^{\infty} \frac{x_j^d}{j^{s_d}}.  
\]

Let $k \in\{1,\ldots, D\}$.  For all 
\[
\Delta = (d_1, d_2, \ldots, d_k) \in \mcd_k  \qqand   J = (j_1, j_2, \ldots, j_k) \in \N^k 
\]
we define the multiplicative coefficient 
\begin{align*}
a_{\Delta,J} 
& = a_{(d_1, d_2, \ldots, d_k),(j_1, j_2, \ldots, j_k)}  = a_{d_1,j_1} a_{d_2,j_2} \cdots a_{d_k,j_k} \\ 
& = \frac{ 1} { j_1^{s_{d_1}} \cdots j_k^{s_{d_k}} } 
\end{align*}.
Then 
\begin{align*} 
 \sum_{J \in \N^k}  a_{\Delta,J} \  x_J^{\Delta} 
& = \left( \mba_{d_1}, \mbx^{d_1 } \right) \left( \mba_{d_2}, \mbx^{d_2 } \right) 
\cdots \left( \mba_{d_k}, \mbx^{d_k } \right) \\
& =  \sum_{(j_1,\ldots, j_k) \in \N^k} 
\frac{ x_{j_1}^{d_1} \cdots x_{j_k}^{d_k} } { j_1^{s_1} \cdots j_k^{s_k} }. 
\end{align*}
We obtain the $q$-multiplicative polynomial of degree at most $D$ 
\begin{align*}
P(\mbx) & = P(\mbx,s_1,\ldots, s_D) \\
& = \sum_{k=1}^D \sum_{ \Delta \in \mcd_k} \sum_{J \in \N^k}  a_{\Delta,J} \  x_J^{\Delta} \\
& =   \sum_{k=1}^D \sum_{ (d_1,\ldots, d_k) \in \mcd_k} 
 \sum_{(j_1,\ldots, j_k) \in \N^k} 
 \frac{ x_{j_1}^{d_1} \cdots x_{j_k}^{d_k}}{ j_1^{s_1} \cdots j_k^{s_k} }. 
\end{align*}
This series converges absolutely for all $\mbx = (x_j)_{j=1}^{\infty}  \in \ell^q$. 

For $D=1$ and $q=2$ and for complex numbers $s$ with $\Re(s) > 1/2$, 
the multiplicative polynomial is the Dirichlet series 
\[
L(\mbx) =   \sum_{n=1}^{\infty} \frac{ x_n}{n^s}.   
\]
This series converges absolutely for all $\mbx = (x_n)_{n=1}^{\infty} \in \ell^2$.  
The corresponding linear equation  in infinitely many variables  is simply 
\[
 \sum_{n=1}^{\infty} \frac{ x_n}{n^s} = b.
 \]
Applying Theorem~\ref{MP:theorem:linear-approx}, 
we obtain the following result, which asserts the existence of a Dirichlet series 
taking prescribed values at a given infinite set of points in the half-plane $\Re(s) > 1/2$.

\bt                  \label{MP:theorem:Dirichlet}
Let $I$ be an infinite set and let $(b_i)_{i\in I}$ and $(s_i)_{i\in I}$ 
be sequences of complex numbers with $\Re(s_i) > 1/2$ for all $i \in I$.
Consider the linear equation in infinitely many variables 
\[
L_i(\mbx) =   \sum_{n=1}^{\infty} \frac{ x_n}{n^{s_i}} = b_i.  
\]
Let $M > 0$. 
If, for every finite subset $S$ of I and every $\varepsilon > 0$, 
the finite set of linear inequalities 
\[
 \left|  \sum_{n=1}^{\infty} \frac{ x_n}{n^{s_i}} -b_i \right| \leq \varepsilon 
 \qquad \text{for all $i \in S$}   
\]
has a solution $\mbx_{S,\varepsilon} \in \ell^2$ 
with $\| \mbx_{S,\varepsilon} \|_2 \leq M$, 
then the infinite set of linear equations 
\[
  \sum_{n=1}^{\infty} \frac{ x_n}{n^{s_i}}  = b_i 
 \qquad \text{for all $i \in I$} 
\] 
has an exact solution $\mbx \in \ell^2$ with $\| \mbx \|_2 \leq M$. 
\et

\section{Open problems}.     \label{MP:OpenProblems}
The results in this paper suggest several questions.

\benum
\item
Are Theorems~\ref{MP:theorem:multiplicative},~\ref{MP:theorem:multiplicative-approx},  and~\ref{MP:theorem:multiplicative-approx-m}  true for  
 infinite sets of polynomial equations in which the 
polynomials are not multiplicative? 

\item
Are there classes \mcc\ of infinite sets of polynomial equations in infinitely many variables 
for which there is an integer $S = S(\mcc)$ such that, if every finite set of at most $S$ equations 
has a solution or an approximate solution, then the infinite set of equations has a solution?

\item
Let $E$ be a subfield of the complex numbers.  
Let $\{P_i(\mbx) = b_i: i \in I\}$ be an infinite set of polynomial equations 
in which every finite subset of the equations has an exact solution 
or an approximate solution in $E$.
 \benum
\item
For what subfields $E$ of the complex numbers is it true that
the infinite set of equations have an exact solution in $E$?  
 \item
For what sets of polynomial equations might this be true?  
\eenum

\item 
The ``finitely many implies infinitely many" paradigm applies to linear equations 
and multiplicative polynomial equations with bounded or sequentially bounded norms 
(Theorems~\ref{MP:theorem:multiplicative},~\ref{MP:theorem:multiplicative-approx},  
and~\ref{MP:theorem:multiplicative-approx-m}).  For equations of this special type, 
countably many implies finitely many.  But there may be other classes of equations 
for which this paradigm does not hold, but a stronger condition (``countably many 
implies uncountably many") is true.  
The problem is to determine if such ``new" classes of equations exist. 
\item
Is there an uncountably infinite set of equations such that every countably infinite subset 
of the equations has an exact or an approximate solution, but the 
uncountably infinite set of equations has no exact solution?  

\eenum

\appendix

\section{The theorems of F. Riesz and Abian-Eslami}         \label{MP:appendix:Riesz}

\bt
Let $(p,q)$ be a conjugate pair and let $M > 0$.  
Let $\mba_i = (a_{i,j})_{j=1}^{\infty} \in \ell^p$  for all $i  \in \N$ 
and let $(b_i)_{i \in \N}$ be a sequence of numbers.   
The following are equivalent.

\benum
\item[(a)]
For all $r \in \N$, 
there exists $\mbx_r  = (x_{r,j})_{j=1}^{\infty} \in \ell^q$ such that 
\[
\|\mbx_r\|_q \leq M
\qqand 
 \sum_{j=1}^{\infty} a_{i,j} x_{r,j} = b_i \qquad \text{for all $i \in \{ 1,2,3,\ldots, r\}$.}
\]
\item[(b)]
For all $r \in \N$ and $h_1,\ldots, h_r \in \R$, 
\[
\left| \sum_{i=1}^r h_i b_i \right| \leq M \left( \sum_{j=1}^{\infty}  
\left| \sum_{i=1}^r h_i  a_{i,j} \right|^p \right)^{1/p}. 
\]
\item[(c)]
There exists $\mbx  = (x_j)_{j=1}^{\infty} \in \ell^q$ such that 
\[
\|\mbx\|_q \leq M
\qqand 
 \sum_{j=1}^{\infty} a_{i,j} x_j  = b_i \qquad \text{for all $i \in \N$.}
\]
\eenum
\et

F. Riesz~\cite{ries13} proved the equivalence of (b) and (c) 
and Abian-Eslami~\cite{abia-esla82} proved (a) implies (c).

\begin{proof}
First we prove that (a) implies (b).
For all  $h_1,\ldots, h_r$, 
\[
  \sum_{i=1}^r h_i  \mba_i =   \sum_{i=1}^r h_i  (a_{i,j})_{j=1}^{\infty} 
  =   \left( \sum_{i=1}^r h_i  a_{i,j} \right)_{j=1}^{\infty} \in \ell^p 
\]
and 
\[
\left\|   \sum_{i=1}^r h_i  \mba_i \right\|_p 
=  \left( \sum_{j=1}^{\infty}  \left| \sum_{i=1}^r h_i  a_{i,j} \right|^p  \right)^{1/p}. 
\]
We have $\mbx_r  = (x_{r,j})_{j=1}^{\infty} \in \ell^q$ and so, 
by H\" older's inequality, 
\[
\left( \sum_{i=1}^r h_i  \mba_i \right) \mbx_r 
= \left( \sum_{i=1}^r h_i  a_{i,j} x_j \right)_{j=1}^{\infty} \in \ell^1
\]
and 
\[
 \left\| \left( \sum_{i=1}^r h_i  \mba_i \right) \mbx_r  \right\|_1 = \sum_{j=1}^{\infty} \left| \sum_{i=1}^r h_i  a_{i,j} x_j \right| < \infty.
\]
For every positive integer $N$, let 
\[
g(N) = \sum_{j=N+1}^{\infty} \left| \sum_{i=1}^r h_i  a_{i,j} x_j \right|.
\]
We have 
\[
\lim_{N\rightarrow\infty} g(N) = 0  
\]
and 
\begin{align*}
\left| \sum_{i=1}^r h_i b_i \right| & = \left| \sum_{i=1}^r h_i \left( \sum_{j=1}^{\infty} a_{i,j} x_{r,j} \right) \right| \\ 
& \leq   
\left| \sum_{i=1}^r h_i   \sum_{j=1}^N a_{i,j} x_{r,j}   \right|
+  \left| \sum_{i=1}^r h_i  \sum_{j=N+1}^{\infty} a_{i,j} x_{r,j}  \right| \\ 
& \leq  
\left|  \sum_{j=1}^N \sum_{i=1}^r h_i   a_{i,j} x_{r,j}   \right| 
+  \sum_{j=N+1}^{\infty}  \left| \sum_{i=1}^r h_i  a_{i,j} x_{r,j}  \right|  \\ 
& \leq   \sum_{j=1}^N  \left| \sum_{i=1}^r h_i   a_{i,j} x_{r,j}   \right| + g(N)  \\
& \leq \sum_{j=1}^{\infty}  \left|  \sum_{i=1}^r h_i   a_{i,j}  x_{r,j}  \right| + g(N)  \\
& = \left\|\left( \sum_{i=1}^r h_i  \mba_i \right) \mbx_r  \right\|_1+ g(N)  \\
& \leq  \left\|\sum_{i=1}^r h_i\mba_i \right\|_p    \|\mbx_r\|_q     +g(N)\\ 
& \leq   M \left( \sum_{j=1}^{\infty}  \left| \sum_{i=1}^r h_i   a_{i,j}   \right|^p\right)^{1/p} +g(N). 
\end{align*}
This inequality is valid for all $N$ and $\lim_{N\rightarrow\infty} g(N) = 0$.
Therefore, 
\[
\left| \sum_{i=1}^r h_i b_i \right| 
\leq M \left( \sum_{j=1}^{\infty}  \left| \sum_{i=1}^r h_i   a_{i,j}   \right|^p\right)^{1/p}.
\]
Thus, (a) implies (b).

 Next we prove that (b) implies (c).
Let $W$ be the vector subspace of $\ell^p$ spanned by the set $\{\mba_i:i \in \N\}$.
Let $h_1,\ldots, h_r,h'_1,\ldots, h'_r$ be scalars such that 
\[
\sum_{i=1}^r h_i\mba_i = \sum_{i=1}^r h'_i\mba_i. 
\]
Then 
\[
\sum_{i=1}^r (h_i  -  h'_i ) \mba_i = \mbo 
\]
and 
\begin{align*}
\left| \sum_{i=1}^r h_i b_i - \sum_{i=1}^r h'_i b_i \right| 
& =  \left| \sum_{i=1}^r (h_i  -  h'_i ) b_i \right| \\ 
& \leq M \left( \sum_{j=1}^{\infty}  \left| \sum_{i=1}^r  (h_i  -  h'_i ) a_{i,j} \right|^p \right)^{1/p}\\ 
& = M   \left\| \sum_{i=1}^r (h_i  -  h'_i ) \mba_i \right\|_p = 0. 
\end{align*}
It follows that there is a well-defined  linear functional $f$ on $W$ such that   
\[
f\left(\sum_{i=1}^r h_i\mba_i \right) = \sum_{i=1}^r h_i b_i.   
\]
In particular, $f_d(\mba_i) = b_i$.  Because 
\[
\left| f\left(\sum_{i=1}^r h_i\mba_i \right) \right| = \left| \sum_{i=1}^r h_i b_i\right| 
\leq M   \left\|\sum_{i=1}^r h_i  \mba_i \right\|_p 
\]
for all $\sum_{i=1}^r h_i  \mba_i  \in W$, the linear functional $f$ has norm $\|f\| \leq M$. 
By the Hahn-Banach theorem, there is a bounded linear functional $F$ on $\ell^p$ 
such that $F(\mbw) = f_d(\mbw)$ for all $\mbw \in W$ and $\| F\| \leq M$.

For every bounded linear functional $F$ on $\ell^p$ there is a sequence $\mbx \in \ell^q$ 
such that  $F(\mba) = (\mba,\mbx)$ for all $\mba \in \ell^p$.  
For all $i \in \N$ we have
\[
b_i = f_d(\mba_i) = F(\mba_i) = (\mba_i,\mbx) = \sum_{j=1}^{\infty} a_{i,j}x_j.
\]
Thus, (b) implies (c).

The proof that (c) implies (a) is immediate.
\end{proof}

\end{document}